%% file: couple.tex
\documentclass[12pt]{article}
\usepackage{amsmath}
\usepackage{amssymb}
\usepackage[dvips]{graphicx}
\usepackage{epsfig}
\usepackage{vector}
\def\eps{\varepsilon}
\font\tencmmib=cmmib10 \skewchar\tencmmib '60
\newfam\cmmibfam
\textfont\cmmibfam=\tencmmib

\def\bbox{\quad\hbox{\vrule \vbox{\hrule \vskip2pt \hbox{\hskip2pt
\vbox{\hsize=1pt}\hskip2pt} \vskip2pt\hrule}\vrule}}
\def\lessim{\ \lower4pt\hbox{$
\buildrel{\displaystyle <}\over\sim$}\ }
\def\gessim{\ \lower4pt\hbox{$\buildrel{\displaystyle >}
\over\sim$}\ }

\def\P{{\cal P}}

\def\A{{\cal A}}

\def\eps{{\varepsilon}}

\def\card{{\mbox{card}}}

\def\la{{\Bigl\langle}}
\def\ra{{\Bigr\rangle}}

\def\qed{\hfill\break\rightline{$\bbox$}}
\newcommand{\e}{\mathbb{E}}

\newcommand{\Reals}{\mathbb{R}}

\newcommand{\vsi}{{\vec{\sigma}}}
\newcommand{\vrho}{{\vec{\rho}}}
\newcommand{\vtau}{{\vec{\tau}}}

\newtheorem{proposition}{Proposition}
\newtheorem{lemma}{Lemma}
\newtheorem{theorem}{Theorem}

\makeatletter
\@addtoreset{equation}{section}

\makeatother

\input invlat.tex

\begin{document}

\title{
A note on the free energy of the coupled system 
in the Sherrington-Kirkpatrick model.
}
\author{Dmitry Panchenko \thanks{
Department of Mathematics, Massachusetts Institute
of Technology, 77 Massachusetts Ave, Cambridge, MA 02139
email: panchenk@math.mit.edu
}\\
{\it Department of Mathematics}\\
{\it Massachusetts Institute of Technology}\\
}

\maketitle
\begin{abstract}
In this paper we consider a system of spins that consists
of two configurations $\vsi^1,\vsi^2\in\Sigma_N=\{-1,+1\}^N$
with Gaussian Hamiltonians $H_N^1(\vsi^1)$ and $H_N^2(\vsi^2)$
correspondingly, 
and these configurations
are coupled on the set where their overlap is fixed
$\{R_{1,2}=N^{-1}\sum_{i=1}^N \sigma_i^1\sigma_i^2 = u_N\}.$
We prove the existence of the thermodynamic limit of the free energy
of this system given that $\lim_{N\to\infty}u_N = u\in[-1,1]$
and give the analogue of the Aizenman-Sims-Starr
variational principle that describes this limit via random overlap structures.
\end{abstract}

\vspace{0.5cm}


\section{Introduction and main results.}

In this paper we will consider a system that consists
of two configurations of spins that are coupled by
fixing their overlap. Our main goal is to prove the
existence of the thermodynamic limit of the free energy of 
this system and to give the characterization of this limit
via random overlap structures in the sense of 
Aizenman-Sims-Starr \cite{ASS}. Let us start
by introducing all necessary notations and definitions.

For any $N\geq 1,$ let us consider a space $\Sigma_N=\{-1,+1\}^N$
and consider two Hamiltonians $H_N^{\ell}(\vsi)$ for $\ell=1,2$
on $\Sigma_N$ given by
\begin{equation}
H_{N}^{\ell}(\vsi)=N^{1/2}\sum_{p\geq 1}\frac{a_p^{\ell}}{N^{p/2}}
\sum_{i_1,\ldots,i_p}g_{i_1,\ldots,i_p} \sigma_{i_1}\ldots\sigma_{i_p},
\label{Ham}
\end{equation}
where $(g_{i_1,\ldots,i_p})$ are standard Gaussian random variables 
independent for all $p\geq 1$ and all $(i_1,\ldots,i_p),$
and the sequences $(a_p^{\ell})_{p\geq 1}$ are such that
\begin{equation}
\sum_{p\geq 1} (a_p^{\ell})^2 <\infty.
\label{L2}
\end{equation}
For $\ell,\ell'\in{1,2},$ let us define the functions
$\xi_{\ell,\ell'}:[-1,1]\to\Reals$ by
\begin{equation}
\xi_{\ell,\ell'}(x)=\sum_{p\geq 1} a_p^{\ell} a_p^{\ell'} x^p
\label{xi}
\end{equation}
so that 
\begin{equation}
\frac{1}{N}\e H_N^{\ell}(\vsi^1) H_N^{\ell'}(\vsi^2)
=\xi_{\ell,\ell'}(R_{1,2}),
\label{Rxi}
\end{equation}
where the overlap 
$$
R_{1,2}
=R(\vsi^1,\vsi^2)
=\frac{1}{N}\sum_{i\leq N}\sigma_i^1\sigma_i^2.
$$
The condition (\ref{L2}) implies that 
the functions $\xi_{\ell,\ell'}$ are well-defined 
and smooth on $[-1,1].$ 
From now on we will also assume that 
the sequences $(a_p^{\ell})$ are such that
the functions $\xi_{\ell,\ell'}$ are convex on $[-1,1].$ 
For example, this holds if $a_p^{\ell}=0$ for $p$ odd
and $a_p^{\ell}\geq 0$ for $p$ even. 
We define the functions,
\begin{equation}
\theta_{\ell,\ell'}(x)=x\xi_{\ell,\ell'}'(x)-\xi_{\ell,\ell'}(x).
\label{theta}
\end{equation}
The convexity of $\xi_{\ell,\ell'}$ implies that for any 
$x,y\in [-1,1]$ we have
\begin{equation}
\xi_{\ell,\ell'}(x) - x\xi_{\ell,\ell'}'(y)+\theta_{\ell,\ell'}(y)\geq 0.
\label{positive}
\end{equation}
Given $u\in[-1,1],$ let us consider a sequence $(u_N)_{N\geq 1}$
such that for each $N$ we have $u_N=k/N$ 
for some integer $-N\leq k\leq N$ and such that 
$\lim_{N\to\infty} u_N=u.$ 
Given the external fields $h_1, h_2\in\Reals,$
we define,
\begin{equation}
F_N(u_N)=\frac{1}{N}\e\log Z_N(u_N),
\label{FN}
\end{equation}
where
\begin{equation}
Z_N(u_N)=
\sum_{R_{1,2}=u_N} 
\exp\Bigl(
\sum_{\ell\leq 2} H_N^{\ell}(\vsi^{\ell}) + \sum_{\ell\leq 2}
h_{\ell}
\sum_{i\leq N} \sigma_i^{\ell}
\Bigr). 
\label{ZN}
\end{equation}
The quantity $F_N(u_N)$ represents the free energy 
of the set of configurations $\{R_{1,2}=u_N\}.$
The main reason that $u_N$ was chosen of the type $k/N$ is
that this set be not empty. 

Our first goal will be to prove the following.
\begin{theorem}\label{Th1}
The limit 
\begin{equation}
\lim_{N\to\infty}F_N(u_N) = \P(u)
\label{limit}
\end{equation}
exists and depends on $u$ but not on the sequence $(u_N).$
\end{theorem}

The main idea in the proof of this Theorem is the interpolation
method of Guerra-Toninelli which was developed by authors 
in \cite{GT} to prove the existence of the thermodynamic limit 
of the free energy of one copy of the system with Hamiltonian $H_N^1(\vsi).$ 
They also extended their method in \cite{GT2} to prove the existence
of the thermodynamic limit in a variety of mean field models.
In fact, a part of the proof of Theorem \ref{Th1} is very
similar to the proof of the main result in \cite{GT2}
which was motivated by the idea of restricting to the set
of configurations with given overlap introduced by Michel Talagrand in
\cite{SG}. However, the situation considered in Theorem \ref{Th1} 
is slightly different, mainly, due to the fact that we consider
the set $\{R_{1,2} = u_N\}$ of configurations with overlap
exactly equal to $u_N$ rather than being in the neighborhood of
$u_N.$ This will require some additional approximation result,
Lemma \ref{Lemma1} below.
We will prove that the sequence $F_N(u_N)$ can be 
approximated by a superadditive sequence over the restricted range of 
indices and apply the following Proposition due to DeBruijn-Erd\"os
\cite{DBE} (see also Theorem 1.9.1 in \cite{Steele}).

\begin{proposition}\label{Prop1}(DeBruijn-Erd\"os)
If the sequence $(a_N)$ of real numbers satisfies the superadditivity
condition 
$$
a_{m+n}\geq a_m+a_n \mbox{ over the restricted range }
\frac{1}{2}n\leq m\leq 2n,
$$
then $\lim_{n\to\infty}a_n/n = \sup a_n/n.$
\end{proposition}

Next, we will characterize the limit $\P(u)$ in (\ref{limit})
via the analogue of Aizenman-Sims-Starr variational principle \cite{ASS}.
This characterization is motivated by the following idea.
In \cite{Tnew} Michel Talagrand proved a certain 
replica symmetry breaking upper bound on $F_N(u_N)$ and 
conjectured that the bound should be precise in the limit,
i.e. should be equal to $\P(u)$ in (\ref{limit}).
He also emphasized that the computation of this limit
is a natural approach to solving the so called chaos problem.
It is interesting to note that the formula conjectured by
Talagrand can be written via Derrida-Ruelle probability cascades 
as in the case of the Parisi formula in the Sherrington-Kirkpatrick 
model. On the other hand, the Parisi formula in the SK model
written via Derrida-Ruelle cascades can be included in a broader
variational principle described in \cite{ASS}.
This connection motivates us to give a variational 
characterization of the limit $\P(u)$ 
in terms of random overlap structures in the sense of
Aizenman-Sims-Starr \cite{ASS}.
We hope that this characterization will provide some insight
into what should be the correct Parisi ansatz for $\P(u)$
and whether the formula conjectured by Talagrand indeed holds.

Given a parameter $\delta>0,$ we define 
{\it the random overlap structure} (ROSt) as 
the following collection of: 

(1) a countable set $\A;$

(2) a sequence $(q_{\alpha,\beta}^{\ell,\ell'})$ for
$\alpha,\beta\in\A, \ell,\ell'\in \{1,2\}$ such that 
\begin{equation}
|q_{\alpha,\beta}^{\ell,\ell'}|\leq 1,\,\,
q_{\alpha,\alpha}^{\ell,\ell}=1 
\,\,\mbox{ and }\,\,
|q_{\alpha,\alpha}^{1,2}-u|\leq \delta;
\label{rost2}
\end{equation}

(3) an arbitrary random sequence $(w_{\alpha})_{\alpha\in\A}$
such that 
\begin{equation}
w_{\alpha}\geq 0 \mbox{ and }
\sum_{\alpha\in\A} w_{\alpha}=1 \mbox{ a.s.;}
\label{rost3}
\end{equation}

(4) Gaussian sequences $(z^1(\alpha),z^2(\alpha))_{\alpha\in\A}$ and 
$(y^1(\alpha),y^2(\alpha))_{\alpha\in\A}$ independent of each other and of
the sequence $(w_{\alpha})_{\alpha\in\A}$ with the following covariance
operators
\begin{equation}
\e z^{\ell}(\alpha) z^{\ell'}(\beta) 
= \xi_{\ell,\ell'}'(q_{\alpha,\beta}^{\ell,\ell'})
\,\,\,\mbox{ and }\,\,\,
\e y^{\ell}(\alpha) y^{\ell'}(\beta) =
 \theta_{\ell,\ell'}(q_{\alpha,\beta}^{\ell,\ell'}).
\label{rost4}
\end{equation}
Let $(z_i^1(\alpha),z_i^2(\alpha))_{\alpha\in\A}$ 
be a sequence of independent copies of
$(z^1(\alpha),z^2(\alpha))_{\alpha\in\A}$ for $i\geq 1.$
We also assume that all random variables here are independent
of the Hamiltonians $H_N^{\ell}(\vsi).$ 
Let us denote such generic collection (1) - (4) as $\Omega_{\delta},$
where we will make the dependence of $\Omega_{\delta}$ on 
the parameter $\delta$ in (\ref{rost2}) explicit. 

One could try to describe conditions on the sequence 
$(q_{\alpha,\beta}^{\ell,\ell'})$
that would guarantee the existence of the Gaussian sequences with the covariance
structure (\ref{rost4}). Instead, we will simply assume that we consider only 
random overlap structures $\Omega_{\delta}$
that such sequences exist. One reason why we are not interested in the general
case is because, as in \cite{ASS}, one particular ROSt will play a special
role in characterization of the limit $\P(u)$ in (\ref{limit})
and it will be constructed explicitly. 
Given a ROSt $\Omega_{\delta},$ let us now consider the quantity
\begin{eqnarray}
G_N(u_N,\Omega_{\delta}) 
&=& 
\frac{1}{N}\e\log\sum_{\alpha\in\A}w_{\alpha} \sum_{R_{1,2}=u_N}
\exp \sum_{\ell\leq 2}\sum_{i\leq N} \sigma_i^{\ell}(z_i^{\ell}(\alpha) 
+ h_{\ell})
\nonumber
\\
&-&
\frac{1}{N}\e\log\sum_{\alpha\in\A}w_{\alpha}
\exp \sqrt{N}\sum_{\ell\leq 2}y^{\ell}(\alpha).
\label{GN}
\end{eqnarray}
The following theorem holds.
\begin{theorem}\label{Th2}
There exists a sequence $(u_N')$ such that 
$\lim_{N\to\infty} u_N' = u$ and such that
the limit in (\ref{limit})
\begin{equation}
\P(u) =\lim_{N\to\infty}
\lim_{\delta\to 0}\inf_{\Omega_{\delta}} G_N(u_N',\Omega_{\delta}).
\label{VAR}
\end{equation}
\end{theorem}

\section{Proof of Theorem \ref{Th1}.}

Given $\eps>0,$ let us consider a set
\begin{equation}
U_{N,\eps}=\bigl\{(\vsi^1,\vsi^2): R_{1,2}\in[u_N-\eps,u_N+\eps]\bigr\}
\label{UNeps}
\end{equation}
and define
\begin{equation}
F_N(U_{N,\eps})=\frac{1}{N}\e\log \sum_{U_{N,\eps}} 
\exp\Bigl(
\sum_{\ell\leq 2} H_N^{\ell}(\vsi^{\ell}) + \sum_{\ell\leq 2}
h_{\ell}
\sum_{i\leq N} \sigma_i^{\ell}
\Bigr). 
\label{FUN}
\end{equation}
In order to utilize the ideas of Guerra and Toninelli in 
\cite{GT} and \cite{GT2},
we first need to prove the following approximation result.

\begin{lemma}\label{Lemma1}
There exists a constant $L$ independent of $N$ such that
for all $\eps\in [0,1]$
\begin{equation}
F_N(U_{N,\eps})\leq F_N(u_N) + L\sqrt{\eps}.
\label{Fcompare}
\end{equation}
\end{lemma}
{\bf Proof.} 
For each $\vsi^1\in\Sigma_N$ let us consider the sets
$$
U_{\eps}(\vsi^1)=\bigl\{\vsi^2: R_{1,2}\in[u_N-\eps,u_N+\eps]\bigr\},
\,\,\,
U(\vsi^1)=
\bigl\{\vsi^2: R_{1,2}=u_N\bigr\}.
$$
For each $\vsi^2\in U_{\eps}(\vsi^1)$ we can find an element
$\pi(\vsi^1,\vsi^2)\in U(\vsi^1)$ such that the Hamming distance
\begin{equation}
d(\vsi^2, \pi(\vsi^1,\vsi^2))=\frac{1}{N}\sum_{i\leq N}
I(\sigma_i^2\not = \pi(\vsi^1,\vsi^2)_i)\leq \frac{\eps}{2}.
\label{pi}
\end{equation}
Indeed, since $R_{1,2}=1-2d(\vsi^1,\vsi^2),$ for 
$\vsi^2\in U_{\eps}(\vsi^1)$ we have
$$
\frac{1-u_N}{2}-\frac{\eps}{2}\leq 
\frac{1}{N}\sum_{i\leq N} I(\sigma_i^1\not =\sigma_i^2) \leq
\frac{1-u_N}{2}+\frac{\eps}{2}.
$$
Therefore, by changing at most $N\eps/2$ coordinates of the vector
$\vsi^2$ we can obtain a vector $\pi(\vsi^1,\vsi^2)$ such that
$$
\frac{1}{N}\sum_{i\leq N} I(\sigma_i^1\not = \pi(\vsi^1,\vsi^2)_i)
= \frac{1-u_N}{2},
$$
which means that $\pi(\vsi^1,\vsi^2)\in U(\vsi^1)$
and $d(\vsi^2,\pi(\vsi^1,\vsi^2))\leq \eps/2.$
If we write
\begin{eqnarray*}
H_N^{2}(\vsi^{2}) + h_{2} \sum_{i\leq N} \sigma_i^{2}
&=&
H_N^{2}(\pi(\vsi^{1},\vsi^2)) + h_2 \sum_{i\leq N} 
\pi(\vsi^1,\vsi^2)_i 
\\
&+&
H_N^{2}(\vsi^{2}) - H_N^{2}(\pi(\vsi^{1},\vsi^2))
+ h_2 \sum_{i\leq N} (\sigma_i^{2} - \pi(\vsi^1,\vsi^2)_i)
\end{eqnarray*}
then, clearly,
$F_N(U_{N,\eps})\leq \mbox{I} + \mbox{II},$
where
\begin{eqnarray*}
\mbox{I} 
&=&
\frac{1}{N}\e\max_{U_{N,\eps}} 
\Bigl(
H_N^{2}(\vsi^{2}) - H_N^{2}(\pi(\vsi^{1},\vsi^2))
+ h_2 \sum_{i\leq N} (\sigma_i^{2} - \pi(\vsi^1,\vsi^2)_i)
\Bigr)
\\
&\leq&
\frac{1}{N}\e\max_{U_{N,\eps}} 
\bigl(
H_N^{2}(\vsi^{2}) - H_N^{2}(\pi(\vsi^{1},\vsi^2))
\bigr)
+|h_2|\eps
\end{eqnarray*}
and
$$
\mbox{II} = 
\frac{1}{N}\e\log \sum_{U_{N,\eps}} \exp
\Bigl(
H_N^{1}(\vsi^{1}) + \sum_{i\leq N} h_{1} \sigma_i^{1}
+
H_N^{2}(\pi(\vsi^{1},\vsi^2)) + h_2 \sum_{i\leq N}  
\pi(\vsi^1,\vsi^2)_i 
\Bigr).
$$
To estimate the first term in I we use 
Slepian's inequality that implies (see \cite{LT})
\begin{eqnarray*}
&&
\e\max_{U_{N,\eps}} 
\bigl(
H_N^{2}(\vsi^{2}) - H_N^{2}(\pi(\vsi^{1},\vsi^2))
\bigr)
\\
&&
\leq
3
\sqrt{\log\card U_{N,\eps}}
\max_{U_{N,\eps}} 
\Bigl(
\e\bigl(
H_N^{2}(\vsi^{2}) - H_N^{2}(\pi(\vsi^{1},\vsi^2))
\bigr)^2 
\Bigr)^{1/2}
\\
&&
\leq
6N\sqrt{\log 2}\max_{U_{N,\eps}} 
\Bigl(
\xi_{2,2}(1)-\xi_{2,2}\bigl( R(\vsi^2,\pi(\vsi^1,\vsi^2)) \bigr)
\Bigr)^{1/2},
\end{eqnarray*}
where we used (\ref{Rxi}) and an estimate 
$\card U_{N,\eps} \leq 2^{2N}.$ By (\ref{pi})
$$
R(\vsi^2,\pi(\vsi^1,\vsi^2))
=1-2d(\vsi^2,\pi(\vsi^1,\vsi^2))\geq 1 - \eps.
$$
Therefore,
$$
\Bigl|\xi_{2,2}(1)-\xi_{2,2}\bigl( R(\vsi^2,\pi(\vsi^1,\vsi^2)) 
\bigr)\Bigr| \leq \max_{x\in[-1,1]}|\xi_{2,2}'(x)|\eps
$$
and, thus, $\mbox{I}\leq L\sqrt{\eps}.$
To estimate II we will simply count how many elements
$\vsi\in U_{\eps}(\vsi^1)$
are projected onto an element $\vsi^2\in U(\vsi^1),$ i.e. 
for $\vsi^2\in U(\vsi^1)$ we consider
$$
\ell(\vsi^1,\vsi^2)
=\card\{\vsi\in U_{\eps}(\vsi^1) : \pi(\vsi^1,\vsi) = \vsi^2\}.
$$
Then, obviously,
\begin{eqnarray*}
\mbox{II} 
&=& 
\frac{1}{N}\e\log \sum_{R_{1,2}=u_N} 
\ell(\vsi^1,\vsi^2)
\exp
\Bigl(
H_N^{1}(\vsi^{1}) + h_1 \sum_{i\leq N}  \sigma_i^{1}
+
H_N^{2}(\vsi^2) + h_2 \sum_{i\leq N}  \sigma_i^2 
\Bigr)
\\
&\leq&
F_N(u_N) + \frac{1}{N}\max_{R_{1,2}=u_N} \log \ell(\vsi^1,\vsi^2).
\end{eqnarray*}
Since by (\ref{pi}), $d(\vsi^2,\pi(\vsi^1,\vsi^2))\leq \eps/2,$
we have
\begin{eqnarray*}
\ell(\vsi^1,\vsi^2)
&\leq& 
\card\{\vsi\in\Sigma_N : 
d(\vsi,\vsi^2)\leq \eps/2\}
=
\card\{\vsi\in\Sigma_N : \sum_{i\leq N} I(\sigma_i\not = 1) 
\leq N\eps/2
\}
\\
&=&
\card\{\vsi\in\Sigma_N : \sum_{i\leq N} \sigma_i
\geq N(1-\eps)\}\leq
2^N \exp(-NI(1-\eps)),
\end{eqnarray*}
where $I(x)=\frac{1}{2}((1+x)\log(1+x)+(1-x)\log(1-x)).$
In the last inequality we used a large deviation estimate
for the Bernoulli r.v. (see, for example, A.9 in \cite{SG}). 
Hence,
$$
\frac{1}{N}\max_{R_{1,2}=u_N} \log \ell(\vsi^1,\vsi^2)
\leq
\log 2 - I(1-\eps)
=
\log\Bigl(1+\frac{\eps}{2-\eps}\Bigr)+\frac{\eps}{2}\log
\frac{2-\eps}{\eps}
\leq L\sqrt{\eps}
$$
for $\eps\in[0,1].$ This finishes the proof of Lemma \ref{Lemma1}.
\qed

Clearly, Lemma \ref{Lemma1} implies that
$$
|F_N(u_N) - F_N(u_N')|\leq L|u_N - u_N'|^{1/2}
$$
for $|u_N - u_N'|\leq 1$ and, therefore, in order to prove
the existence of the limit $\lim_{N\to\infty} F_N(u_N)$ 
for any sequence $(u_N)$ such that $\lim_{N\to\infty}=u$
it is enough to prove it for one such sequence. 
Therefore, from now on we will make a specific choice of
$(u_N)$ that satisfies the following condition,
\begin{equation}
|u_N - u|\leq \frac{1}{N}.
\label{near}
\end{equation}
Clearly, it is possible to take $u_N$ of the type
$u_N=k/N$ that satisfies this condition.  

The next Lemma is similar to the techniques in \cite{GT2}. 

\begin{lemma}\label{Lemma2}
If $(u_N)$ satisfies (\ref{near}) then there exists a constant $A$
independent of $N$ such that the sequence
$$
a_N = NF_N(u_N) - A N^{1/2}
$$
satisfies superadditivity condition
$$
a_{M+N}\geq a_M + a_N \mbox{ over the restricted range }
\frac{1}{2}N\leq M\leq 2N.
$$
\end{lemma}
{\bf Proof.}
Given $N,M\geq 1,$ let us consider a space $\Sigma_{M+N}$
and for each $\vsi\in\Sigma_{N+M}$ we will write
$\vsi=(\vrho,\vtau)$ where 
\begin{equation}
\vrho=(\rho_1,\ldots,\rho_M)=(\sigma_1,\ldots,\sigma_M)\in\Sigma_M,\,\,\,
\vtau=(\tau_1,\ldots,\tau_N)=(\sigma_{M+1},\ldots,\sigma_{M+N})\in\Sigma_N.
\label{srt}
\end{equation}
For $\vsi^1=(\vrho^1,\vtau^1)$ and $\vsi^2=(\vrho^2,\vtau^2)$ we define
$$
R_{1,2}^1=R(\vrho^1,\vrho^2)=\frac{1}{M}\sum_{i\leq M}\rho_i^1\rho_i^2
\,\,\,\mbox{ and }\,\,\,
R_{1,2}^2=R(\vtau^1,\vtau^2)=\frac{1}{N}\sum_{i\leq N}\tau_i^1\tau_i^2.
$$
Let us write the overlap $R_{1,2}=R(\vsi^1,\vsi^2)$ as
$$
R_{1,2}=\frac{M}{M+N}R_{1,2}^1 + \frac{N}{M+N}R_{1,2}^2.
$$
Then we have,
\begin{equation}
U_{M,N}:=\{R_{1,2}^1=u_M, R_{1,2}^2= u_N\}\subseteq
\Bigl\{R_{1,2}= 
u_{M+N}'
:=
\frac{M}{M+N} u_M + \frac{N}{M+N} u_N  
\Bigr\}.
\label{overlapconv}
\end{equation}
If we define $\eps=|u_{M+N}' - u_{M+N}|$ then 
$$
\{R_{1,2}=u_{M+N}'\}\subseteq
U_{M+N,\eps}=\Bigl\{ R_{1,2}\in
[u_{M+N}-\eps,u_{M+N}+\eps]\Bigr\}
$$
and, therefore, $U_{M,N}\subseteq U_{M+N,\eps}.$
This together with Lemma \ref{Lemma1} implies,
$$
F_{M+N}(u_{M+N})\geq F_{M+N}(U_{M+N,\eps}) - L\sqrt{\eps}\geq
F_{M+N}(U_{M,N}) - L\sqrt{\eps},
$$
where
$$
F_{M+N}(U_{M,N})= \frac{1}{M+N}\e\log \sum_{U_{M,N}} 
\exp\Bigl(
\sum_{\ell\leq 2} H_{M+N}^{\ell}(\vsi^{\ell}) + \sum_{\ell\leq 2}
h_{\ell}
\sum_{i\leq M+N} \sigma_i^{\ell}
\Bigr). 
$$
Condition (\ref{near}) implies that $\eps\leq 3/(M+N)$
and, therefore,
\begin{equation}
F_{M+N}(u_{M+N})\geq F_{M+N}(U_{M,N}) - \frac{L}{(M+N)^{1/2}}.
\label{Fstep1}
\end{equation}
Given $t\in[0,1],$ let us consider an interpolating Hamiltonian
$$
H_t(\vsi^1,\vsi^2) = 
\sqrt{t}\sum_{\ell\leq 2} H_{M+N}^{\ell}(\vsi^{\ell}) 
+
\sqrt{1-t}
\sum_{\ell\leq 2} 
\Bigl(
H_{M}^{\ell}(\vrho^{\ell})
+
H_{N}^{\ell}(\vtau^{\ell})
\Bigr)
$$
where the Hamiltonians $H_M^{\ell}, H_N^{\ell}$
and $H_{M+N}^{\ell}$ are independent of each other,
and define a function $\varphi(t)$ by
\begin{equation}
(M+N)\varphi(t)=\e\log \sum_{U_{M,N}} 
\exp\Bigl(
H_t(\vsi^1,\vsi^2)
+ \sum_{\ell\leq 2}
h_{\ell}
\sum_{i\leq M+N} \sigma_i^{\ell}
\Bigr).
\end{equation}
It is easy to see that 
$$
\varphi(1)=F_{M+N}(U_{M,N})
\mbox{ and }
\varphi(0)=\frac{M}{M+N}F_M(u_M) + \frac{N}{M+N}F_N(u_N).
$$
We will show below that for some constant $L,$
\begin{equation}
\varphi'(t)\geq -\frac{L}{N+M}.
\label{Fstep2}
\end{equation}
This control of the derivative will imply that
$\varphi(1)\geq \varphi(0)-L/(M+N)$
and, combining this with (\ref{Fstep1}), we get
$$
(M+N)F_{M+N}(u_{M+N})\geq M F_M(u_M) + N F_N(u_N) - L(M+N)^{1/2}.
$$
If given $A>0$ we consider a sequence $a_N=NF_N(u_N)-AN^{1/2}$ then
this can be written equivalently as,
$$
a_{M+N}\geq a_{M} + a_{N} + A M^{1/2}+ A N^{1/2}-(A+L)(M+N)^{1/2}.
$$
When $N/2\leq M\leq 2N,$ we have
$$
M^{1/2}+N^{1/2} \geq 
(M+N)^{1/2}\Bigl(\sqrt{\frac{1}{3}}+\sqrt{\frac{2}{3}}\Bigr)
$$
and, thus,
$$
A M^{1/2}+ A N^{1/2}-(A+L)(M+N)^{1/2}\geq
\Bigl(\Bigl(\sqrt{\frac{1}{3}}
+\sqrt{\frac{2}{3}}-1\Bigr)A-L\Bigr)(M+N)^{1/2}\geq 0,
$$
if $A$ is large enough. This proves that 
$$
a_{M+N}\geq a_M + a_N 
\mbox{ over the restricted range }
\frac{1}{2}N\leq M\leq 2N,
$$
which is precisely the statement of Lemma. I
Hence, it remains to prove (\ref{Fstep2}).

Let us denote by $\la \cdot \ra_t$ the average with respect to
the Gibbs' measure $G_{M,N}$ on $U_{M,N}$ with Hamiltonian
$$
H_t(\vsi^1,\vsi^2)
+ \sum_{\ell\leq 2}
h_{\ell}
\sum_{i\leq M+N} \sigma_i^{\ell}.
$$
Then the standard computation utilizing Gaussian integration by parts
gives (see, for example, \cite{GT} or Theorem 2.10.1 in \cite{SG}),
\begin{eqnarray*}
&&
(M+N)\varphi'(t)=
\frac{1}{2}\sum_{\ell,\ell'\leq 2}
\e\Bigl\la 
(M+N)\xi_{\ell,\ell'}(R_{\ell,\ell'})-
M\xi_{\ell,\ell'}(R_{\ell,\ell'}^1)-
N\xi_{\ell,\ell'}(R_{\ell,\ell'}^2)
\Bigr\ra_t
\\
&&
-
\frac{1}{2}\sum_{\ell,\ell'\leq 2}
\e\Bigl\la 
(M+N)\xi_{\ell,\ell'}(R(\vsi^{\ell},\bar{\vsi}^{\ell'}))-
M\xi_{\ell,\ell'}(R(\vrho^{\ell},\bar{\vrho}^{\ell'}))-
N\xi_{\ell,\ell'}(R(\vtau^{\ell},\bar{\vtau}^{\ell'}))
\Bigr\ra_t \, ,
\end{eqnarray*}
where $(\bar{\vsi}^1,\bar{\vsi}^2)$ is an independent copy 
of $(\vsi^1,\vsi^2)$ with respect to the Gibbs' measure $G_{M,N}$.
Since
$$
R(\vsi^{\ell},\bar{\vsi}^{\ell'})= 
\frac{M}{M+N}R(\vrho^{\ell},\bar{\vrho}^{\ell'})
+\frac{N}{M+N}R(\vtau^{\ell},\bar{\vtau}^{\ell'}),
$$
the convexity of $\xi_{\ell,\ell'}$ implies that
$$
(M+N)\xi_{\ell,\ell'}(R(\vsi^{\ell},\bar{\vsi}^{\ell'}))
\leq
M\xi_{\ell,\ell'}(R(\vrho^{\ell},\bar{\vrho}^{\ell'}))
+
N\xi_{\ell,\ell'}(R(\vtau^{\ell},\bar{\vtau}^{\ell'})),
$$
and, therefore,
$$
(M+N)\varphi'(t)\geq
\frac{1}{2}\sum_{\ell,\ell'\leq 2}
\e\Bigl\la 
(M+N)\xi_{\ell,\ell'}(R_{\ell,\ell'})-
M\xi_{\ell,\ell'}(R_{\ell,\ell'}^1)-
N\xi_{\ell,\ell'}(R_{\ell,\ell'}^2)
\Bigr\ra_t \, .
$$
For $\ell=\ell'$ we have 
$R_{\ell,\ell}=R_{\ell,\ell}^1=R_{\ell,\ell}^2=1.$
Also since the average $\la\cdot\ra_t$ is defined on $U_{M,N}$
we have $R_{1,2}^1=u_M,$ $R_{1,2}^2 = u_N$ and 
by (\ref{overlapconv}) $R_{1,2}=u_{M+N}'.$
Thus,
$$
(M+N)\varphi'(t)\geq
(M+N)\xi_{1,2}(u_{M+N}') - M\xi_{1,2}(u_M) - N\xi_{1,2}(u_N).
$$
Condition (\ref{near}) implies that for all $N$ we have
$|\xi_{1,2}(u_N)-\xi_{1,2}(u)|\leq L/N$
and this, clearly, implies (\ref{Fstep2}).
\qed

Combining Lemma \ref{Lemma2} and Proposition \ref{Prop1} 
proves that the limit
$\lim_{N\to\infty} a_N/N$ exists and it is, obviously,
equal to the limit $\lim_{N\to\infty} F_N(u_N),$  
which finishes the proof of Theorem \ref{Th1}.

\section{Proof of Theorem \ref{Th2}.}

In this section we will assume that the sequence
$(u_N)$ satisfies (\ref{near}).
Let us start by proving the following upper bound.

\begin{lemma}\label{Lemma3}
For some constant $L$ independent of $N$ we have,
\begin{equation}
F_N(u_N)\leq \inf_{\Omega_{\delta}}G_N(u_N,\Omega_{\delta})
+ L\delta + L N^{-1} .
\label{upper}
\end{equation}
\end{lemma}
{\bf Proof.} 
Consider an arbitrary random overlap structure $\Omega_{\delta}.$
Given $t\in[0,1],$ let us consider a Hamiltonian
$H_t(\alpha,\vsi^1,\vsi^2)$ on the set $\A\times\{R_{1,2}=u_N\}$
given by
\begin{eqnarray*}
H_t(\alpha,\vsi^1,\vsi^2)
&=&
\sqrt{t}\Bigl(
\sum_{\ell\leq 2}H_N^{\ell}(\vsi^{\ell})+ 
\sqrt{N}\sum_{\ell\leq 2}y^{\ell}(\alpha)
\Bigr)
\\
&&
+\sqrt{1-t}\sum_{\ell\leq 2}\sum_{i\leq N}
\sigma_i^{\ell} z_i^{\ell}(\alpha) + 
\sum_{\ell\leq 2}h_{\ell}\sum_{i\leq N}\sigma_i^{\ell},
\end{eqnarray*}
and consider a function
$$
\varphi(t)=\frac{1}{N}\e\log\sum_{\alpha\in\A}w_{\alpha}
\sum_{R_{1,2}=u_N}\exp H_t(\alpha,\vsi^1,\vsi^2).
$$
Clearly, the statement of lemma is then equivalent to
$$
\varphi(1)\leq \varphi(0) + L \delta + L N^{-1}.
$$ 
We will prove this by showing that the derivative
$\varphi'(t)\leq L\delta + LN^{-1}.$
Let us denote by $\la\cdot\ra_t$ the average with respect to 
the Gibbs' measure on $\A\times\{R_{1,2}=u_N\}$ with Hamiltonian
$H_t(\alpha,\vsi^1,\vsi^2).$
Then the standard computation utilizing Gaussian integration by parts
and covariance structure (\ref{rost4}) gives
\begin{eqnarray}
\varphi'(t)
&=&
\frac{1}{2}\sum_{\ell,\ell'\leq 2}
\e\Bigl\la 
\xi_{\ell,\ell'}(R_{\ell,\ell'})-
R_{\ell,\ell'} \xi_{\ell,\ell'}'(q_{\alpha,\alpha}^{\ell,\ell'})
+\theta_{\ell,\ell'}(q_{\alpha,\alpha}^{\ell,\ell'})
\Bigr\ra_t
\label{twoterms}
\\
&-&
\frac{1}{2}\sum_{\ell,\ell'\leq 2}
\e\Bigl\la 
\xi_{\ell,\ell'}(R(\vsi^{\ell},\bar{\vsi}^{\ell'}))-
R(\vsi^{\ell},\bar{\vsi}^{\ell'}) 
\xi_{\ell,\ell'}'(q_{\alpha,\beta}^{\ell,\ell'})
+\theta_{\ell,\ell'}(q_{\alpha,\beta}^{\ell,\ell'})
\Bigr\ra_t \, ,
\nonumber
\end{eqnarray}
where $(\beta,\bar{\vsi}^1,\bar{\vsi}^2)$ is an independent
copy of $(\alpha,\vsi^1,\vsi^2).$
Using the fact that
the average $\la\cdot\ra_t$ is taken over the set where $R_{1,2}=u_N,$
the first sum on the right hand side is equal to 
$$
\e\Bigl\la\xi_{1,2}(u_N)-u_N\xi_{1,2}'(q_{\alpha,\alpha}^{1,2}) 
+ \theta_{1,2}(q_{\alpha,\alpha}^{1,2}) 
\Bigr\ra_t \leq L\delta + LN^{-1},
$$ 
where the last inequality follows from the fact that 
by (\ref{rost2}) we have $|q_{\alpha,\alpha}^{1,2} - u| \leq \delta$ 
and by (\ref{near}) we have $|u_N -u|\leq N^{-1}.$
The second line in (\ref{twoterms})
is negative by (\ref{positive}) and this finishes the proof. 
\qed

To prove the lower bound, let us start with a couple of simple lemmas.

\begin{lemma}\label{Lemma4}
If a sequence $(a_N)$ is such that $\lim_{N\to\infty}a_N/N = \gamma$
then for any $N\geq 1$ we have
$$
\frac{1}{N}\liminf_{M\to\infty}(a_{M+N}-a_M)\leq \gamma.
$$
\end{lemma}
{\bf Proof.} Suppose that for some $N\geq 1$ and for some $\eps>0$
$$
\frac{1}{N}\liminf_{M\to\infty}(a_{M+N}-a_M)\geq \gamma + \eps.
$$
Then there exists $M_0\geq 1$ such that for all $M\geq M_0$
$$
\frac{1}{N}(a_{M+N}-a_M)\geq \gamma + \frac{\eps}{2}
$$
and, therefore, for $k\geq 0$
$$
\frac{1}{N}(a_{M+(k+1)N}-a_{M+kN})\geq \gamma + \frac{\eps}{2}.
$$
Adding these inequalities for $0\leq k\leq m-1$ we get
$$
\frac{1}{N}(a_{M+mN}-a_{M})\geq m\Bigl(\gamma + \frac{\eps}{2}\Bigr)
\,\,\,\mbox{ or }\,\,\,
\frac{1}{mN}(a_{M+mN}-a_{M})\geq \gamma + \frac{\eps}{2}.
$$
Letting $m\to\infty$ yields that 
$\liminf_{N\to\infty} a_N/N \geq \gamma+\eps/2$ and
this contradicts the fact that $\lim_{N\to\infty}a_N/N = \gamma.$
\qed

\begin{lemma}\label{Lemma5}
Consider a sequence $(u_N)$ such that (\ref{near}) holds.
Then there exists a sequence $(u_N')$ such that 
$|u_N' - u|\leq 2/N$ and such that for each $N\geq 1,$
\begin{equation}
\frac{M}{M+N}u_{M} + \frac{N}{M+N} u_N' = u_{M+N}
\label{uprime}
\end{equation}
for infinitely many $M\geq 1.$ 
\end{lemma}
{\bf Proof.}
For a fixed $N,$ consider a sequence $u_N'(M)$ defined
by (\ref{uprime}), i.e.
$$
Nu_N'(M) = (M+N)u_{M+N} - Mu_M.
$$
We have
$$
N(u_N'(M) - u) = (M+N)(u_{M+N} - u) - M(u_M - u)
$$
and, therefore, (\ref{near}) implies that  $N|u_N'(M) - u| \leq 2.$
Since $Nu_N'(M)$ is an integer between $-N$ and $N,$ 
it can take a finite number of values and, thus, we can
find an infinite subsequence $(M_k)_{k\geq 1}$ such that 
$u_N'(M_k)=u_{N}'(M_1).$ Take $u_N' = u_N'(M_1).$
\qed

\begin{theorem}\label{Th3}
There exists a sequence $(u_N')$ such that
$|u_N' - u|\leq L/N$ and for all $N\geq 1,$
$$
\P(u)\geq \lim_{\delta\to 0}\inf_{\Omega_{\delta}}
G_N(u_N',\Omega_{\delta}).
$$ 
\end{theorem}
{\bf Proof.}
If we consider a sequence $a_N = NF_N(u_N)$ then, by Theorem \ref{Th1},
we have that
the limit $\lim_{N\to\infty} a_N/N = \P(u).$ Lemma \ref{Lemma4} then
implies that for any $N\geq 1,$
\begin{equation}
\frac{1}{N}\liminf_{M\to\infty}
\Bigl((M+N)F_{M+N}(u_{M+N}) - MF_M(u_M)\Bigr)\leq \P(u).
\label{ONE}
\end{equation}
We can write
\begin{equation}
\frac{1}{N}\Bigl((M+N)F_{M+N}(u_{M+N}) - MF_M(u_M)\Bigr)=
\frac{1}{N}\e\log Z_{M+N}(u_{M+N}) - \frac{1}{N}\e\log Z_{M}(u_M),
\label{TWO}
\end{equation}
where $Z_{N}(u_N)$ was defined in (\ref{ZN}).
For $\vsi\in\Sigma_{M+N}$ we will write $\vsi=(\vrho,\vtau)$
as in (\ref{srt}).
Consider the sequence $(u_N')$ as in Lemma \ref{Lemma5}.
Then as in (\ref{overlapconv}) the condition (\ref{uprime})
implies that for infinitely many $M\geq 1$ we have
$$
\{R_{1,2}=u_{M+N}\}\supseteq U_{M,N}' :=
\{R_{1,2}^1 = u_M, R_{1,2}^2 = u_N'\}.
$$
For simplicity of notations let us assume that this holds for
all $M\geq 1$ rather than a subsequence $(M_k).$
Therefore,
$$
Z_{M+N}(u_{M+N})\geq Z_{M,N}(u_M,u_N):=
\sum_{U_{M,N}'} 
\exp\Bigl(
\sum_{\ell\leq 2} H_{M+N}^{\ell}(\vsi^{\ell}) + \sum_{\ell\leq 2}
h_{\ell}\sum_{i\leq M+N} \sigma_i^{\ell}
\Bigr). 
$$
Let us decompose the Hamiltonian in $Z_{M,N}(u_{M},u_N)$ as,
\begin{eqnarray}
\sum_{\ell\leq 2} H_{M+N}^{\ell}(\vsi^{\ell}) + \sum_{\ell\leq 2}
h_{\ell}\sum_{i\leq M+N} \sigma_i^{\ell}
&=&
\sum_{\ell\leq 2} H_{M+N}^{\ell}(\vrho^{\ell}) + \sum_{\ell\leq 2}
h_{\ell}\sum_{i\leq M} \rho_i^{\ell}
\label{decompose}
\\
&+& 
\sum_{\ell\leq 2}\sum_{i\leq N}\tau_i^{\ell}
(Z_i^{\ell}(\vrho^{\ell}) + h_{\ell})
+R(\vsi^1,\vsi^2).
\nonumber
\end{eqnarray}
The first two terms on the right hand side represent 
the part of the Hamiltonian that depends on the first $M$ coordinates
$\vrho$ only, i.e. here
\begin{equation}
H_{M+N}^{\ell}(\vrho^{\ell})=(M+N)^{1/2}
\sum_{p\geq 1}\frac{a_p^{\ell}}{(M+N)^{p/2}}
\sum_{i_1,\ldots,i_p\leq M}g_{i_1,\ldots,i_p} 
\rho_{i_1}^{\ell}\ldots\rho_{i_p}^{\ell}.
\label{HMN}
\end{equation}
The third term consists of the terms in the Hamiltonian 
that depend only on one
of the last $N$ coordinates 
$(\tau_1^{\ell},\ldots,\tau_N^{\ell})$ of $\vsi^{\ell},$  
i.e. 
$$
Z_i^{\ell}(\vrho^{\ell}) = 
(M+N)^{1/2}\sum_{p\geq 1}\frac{a_p^{\ell}}{(M+N)^{p/2}}
\sum_{i_1,\ldots,i_{p-1}\leq M}g_{i_1,\ldots,i_{p-1}}^{(i)} \rho_{i_1}^{\ell}
\ldots\rho_{i_{p-1}}^{\ell},
$$
where
$$
g_{i_1,\ldots,i_{p-1}}^{(i)} =
g_{i,i_1,\ldots,i_{p-1}}+g_{i_1,i,\ldots,i_{p-1}}+\ldots 
+g_{i_1,\ldots,i_{p-1},i}.
$$
Finally, the last term $R(\vsi^1,\vsi^2)$ consists of all the terms
of the Hamiltonian that depend on at least two coordinates in 
$\vtau^{\ell}.$ Note that $R(\vsi^1,\vsi^2)$ is independent of all other terms
in (\ref{decompose}) and, therefore, H\"older's inequality implies that
\begin{eqnarray}
\frac{1}{N}\e\log Z_{M,N}(u_{M},u_N) 
&\geq&
\frac{1}{N}\e \log \sum_{U_{M,N}'} \exp
\Bigl(
\sum_{\ell\leq 2} H_{M+N}^{\ell}(\vrho^{\ell}) 
+
\sum_{\ell\leq 2}
h_{\ell}\sum_{i\leq M} \rho_i^{\ell}
\nonumber
\\
&&
+ 
\sum_{\ell\leq 2}\sum_{i\leq N}\tau_i^{\ell}
(Z_i^{\ell}(\vrho^{\ell}) + h_{\ell})
\Bigr).
\label{tempZ}
\end{eqnarray}
For each $(\vrho^1,\vrho^2)$ let us denote
$$
W(\vrho^1,\vrho^2)=\exp
\Bigl(
\sum_{\ell\leq 2} H_{M+N}^{\ell}(\vrho^{\ell}) 
+ \sum_{\ell\leq 2}
h_{\ell}\sum_{i\leq M} \rho_i^{\ell}
\Bigr)
$$
so that (\ref{tempZ}) becomes
\begin{eqnarray}
\frac{1}{N}\e\log Z_{M,N}(u_{M},u_N) 
\geq
\frac{1}{N}\e \log \sum_{U_{M,N}'} W(\vrho^1,\vrho^2)
\exp
\sum_{\ell\leq 2}\sum_{i\leq N}\tau_i^{\ell}
(Z_i^{\ell}(\vrho^{\ell}) + h_{\ell})
.
\label{tempZ2}
\end{eqnarray}
The sequences $(Z_i^{\ell}(\vrho^{\ell}))$ are independent
for different indices $i,$ and
the covariance operator of $(Z_i^{\ell}(\vrho^{\ell}))$ is given by
$$
\e Z_i^{\ell}(\vrho^{\ell}) Z_i^{\ell'}(\vrho^{\ell'}) = 
\sum_{p\geq 1} \Bigl(
\frac{M}{M+N}\Bigr)^{p-1}
a_p^{\ell}a_p^{\ell'} p (R(\vrho^{\ell},\vrho^{\ell'}))^{p-1}
= \xi_{\ell,\ell'}'(R(\vrho^{\ell},\vrho^{\ell'})) + o_M(1)
$$
as $M\to\infty,$
uniformly over $R(\vrho^{\ell},\vrho^{\ell'})\in [-1,1].$
Therefore, one can substitute (up to a small error) 
the random variables $Z_{i}^{\ell}(\vrho^{\ell})$
in (\ref{tempZ2}) with the random variables
\begin{equation}
z_i^{\ell}(\vrho^{\ell}) = M^{1/2}\sum_{p\geq 1}\frac{a_p^{\ell}}{M^{p/2}}
\sum_{i_1,\ldots,i_{p-1}\leq M}g_{i_1,\ldots,i_{p-1}}^{(i)} \rho_{i_1}^{\ell}
\ldots\rho_{i_{p-1}}^{\ell}
\label{zee}
\end{equation}
with covariance operator
\begin{equation}
\e z_i^{\ell}(\vrho^{\ell}) z_i^{\ell'}(\vrho^{\ell'}) 
= \xi_{\ell,\ell'}'(R(\vrho^{\ell},\vrho^{\ell'})).
\label{covzee}
\end{equation}
Namely, we have,
\begin{eqnarray*}
&&
\frac{1}{N}\e \log \sum_{U_{M,N}'} W(\vrho^1,\vrho^2)\exp
\sum_{\ell\leq 2}
\sum_{i\leq N}\tau_i^{\ell}(Z_i^{\ell}(\vrho^{\ell}) + h_{\ell})
\\
&&
=\frac{1}{N}
\e \log \sum_{U_{M,N}'} W(\vrho^1,\vrho^2)\exp 
\sum_{\ell\leq 2}
\sum_{i\leq N}\tau_i^{\ell}(z_i^{\ell}(\vrho^{\ell}) + h_{\ell})
+ o_M(1),
\end{eqnarray*}
when $M\to\infty.$
This is easy to show by interpolating between $Z_i^{\ell}$ and
$z_i^{\ell}$ via
$$
z_i^{\ell}(\vrho^{\ell},t) = \sum_{p\geq 1}a_p^{\ell}
\Bigl(
\frac{t}{M^{(p-1)/2}}
+\frac{1-t}{(M+N)^{(p-1)/2}}
\Bigr)
\sum_{i_1,\ldots,i_{p-1}\leq M}g_{i_1,\ldots,i_{p-1}}^{(i)} \rho_{i_1}^{\ell}
\ldots\rho_{i_{p-1}}^{\ell}
$$
and considering
$$
\varphi(t)= \frac{1}{N}
\e \log \sum_{U_{M,N}'} W(\vrho^1,\vrho^2)\exp 
\sum_{\ell\leq 2}
\sum_{i\leq N}\tau_i^{\ell}(z_i^{\ell}(\vrho^{\ell},t) + h_{\ell}).
$$
Then it is a straightforward calculation to show that
$\varphi'(t)=o_M(1)$ uniformly for $t\in[0,1].$
Thus, we finally get,
\begin{eqnarray}
&&
\frac{1}{N}\e\log Z_{M,N}(u_{M},u_N) 
\geq
\frac{1}{N}\e \log \sum_{U_{M,N}'} W(\vrho^1,\vrho^2)
\exp
\sum_{\ell\leq 2}
\sum_{i\leq N}\tau_i^{\ell}(z_i^{\ell}(\vrho^{\ell}) + h_{\ell})
-o_M(1)
\nonumber
\\
&&
= \frac{1}{N}
\e \log \sum_{R_{1,2}^1=u_M} W(\vrho^1,\vrho^2)
\sum_{R_{1,2}^2 = u_N'}\exp
\sum_{\ell\leq 2}
\sum_{i\leq N}\tau_i^{\ell}(z_i^{\ell}(\vrho^{\ell}) + h_{\ell})
-o_M(1)
\label{onepart}
\end{eqnarray}
where $z_i^{\ell}(\vrho^{\ell})$ are defined in (\ref{zee}).
Next, let us consider
$$
Z_M(u_M)= 
\sum_{R_{1,2}^1=u_M} 
\exp\Bigl(
\sum_{\ell\leq 2} H_{M}^{\ell}(\vrho^{\ell}) + \sum_{\ell\leq 2}
h_{\ell}\sum_{i\leq M} \rho_i^{\ell}
\Bigr). 
$$
Comparing
$$
H_{M}^{\ell}(\vrho^{\ell})=M^{1/2}\sum_{p\geq 1}\frac{a_p^{\ell}}{M^{p/2}}
\sum_{i_1,\ldots,i_p\leq M}
g_{i_1,\ldots,i_p} \rho_{i_1}^{\ell}\ldots\rho_{i_p}^{\ell},
$$
with $H_{M+N}^{\ell}(\vrho^{\ell})$ in (\ref{HMN}) we can write
\begin{equation}
H_M^{\ell}(\vrho^{\ell}) \stackrel{D}{=}
H_{M+N}^{\ell}(\vrho^{\ell}) + \sqrt{N} Y^{\ell}(\vrho^{\ell})
\label{inD}
\end{equation}
where
$$
Y^{\ell}(\vrho^{\ell}) = \frac{1}{\sqrt{N}} 
\sum_{p\geq 1} \Bigl(
\frac{1}{M^{p-1}}-\frac{1}{(M+N)^{p-1}}
\Bigr)^{1/2} a_p^{\ell} \sum_{i_1,\ldots,i_p\leq M}
\tilde{g}_{i_1,\ldots,i_p} \rho_{i_1}^{\ell}\ldots \rho_{i_p}^{\ell}, 
$$
where $(\tilde{g}_{i_1,\ldots,i_p})$ are i.i.d. Gaussian r.v
independent of the Hamiltonians $H_{M+N}^{\ell}(\vrho^{\ell}).$
Using (\ref{inD}), we can write
\begin{eqnarray}
\frac{1}{N}\e\log Z_M(u_M) 
&=& 
\frac{1}{N}\e \log \sum_{R_{1,2}^{1}=u_M} \exp 
\Bigl(
\sum_{\ell\leq 2}H_{M+N}^{\ell}(\vrho^{\ell}) + 
\sum_{\ell\leq 2}h_{\ell}\sum_{i\leq M}\rho_i^{\ell}
+\sqrt{N} \sum_{\ell\leq 2}Y^{\ell}(\vrho^{\ell})
\Bigr)
\nonumber
\\
&=&
\frac{1}{N}\e \log \sum_{R_{1,2}^{1}=u_M} 
W(\vrho^1,\vrho^2)
\exp \sqrt{N} \sum_{\ell\leq 2}Y^{\ell}(\vrho^{\ell}).
\label{inD2}
\end{eqnarray}
It is easy to compute that the covariance
operator of $(Y^{\ell}(\vrho^{\ell}))$ satisfies
$$
\e Y^{\ell}(\vrho^{\ell}) Y^{\ell'}(\vrho^{\ell'})
= \theta_{\ell,\ell'}(R(\vrho^{\ell},\vrho^{\ell'})) + o_M(1)
$$
as $M\to \infty.$ 
Therefore, as above one can substitute 
(up to a small error) the random variables
$(Y^{\ell}(\vrho^{\ell}))$ in (\ref{inD2}) with the random variables
\begin{equation}
y^{\ell}(\vrho^{\ell}) = 
\sum_{p\geq 1}(p-1)
a_p^{\ell} \sum_{i_1,\ldots,i_p\leq M}
\tilde{g}_{i_1,\ldots,i_p} \rho_{i_1}^{\ell}\ldots \rho_{i_p}^{\ell}, 
\label{ys}
\end{equation}
with covariance operator 
\begin{equation}
\e y^{\ell}(\vrho^{\ell}) y^{\ell'}(\vrho^{\ell'})
= \theta_{\ell,\ell'}(R(\vrho^{\ell},\vrho^{\ell'})).
\label{covys}
\end{equation}
(\ref{inD2}) then gives,
\begin{equation}
\frac{1}{N}\e\log Z_M(u_M) 
=
\frac{1}{N}\e \log \sum_{R_{1,2}^{1}=u_M} 
W(\vrho^1,\vrho^2)
\exp \sqrt{N} \sum_{\ell\leq 2}y^{\ell}(\vrho^{\ell})
+ o_M(1)
\label{part2}
\end{equation}
as $M\to\infty.$
Plugging (\ref{onepart}) and (\ref{part2}) into (\ref{TWO}) and (\ref{ONE}) we
get 
\begin{equation}
\P(u)\geq \liminf_{M\to\infty} G_{M,N}
\label{Plim}
\end{equation}
where
\begin{eqnarray}
G_{M,N}
&=&
\frac{1}{N}
\e \log \sum_{R_{1,2}^1=u_M} W(\vrho^1,\vrho^2)
\sum_{R_{1,2}^2 = u_N'}\exp
\sum_{\ell\leq 2}\sum_{i\leq N}
\tau_i^{\ell}(z_i^{\ell}(\vrho^{\ell}) + h_{\ell})
\nonumber
\\
&&
-
\frac{1}{N}\e \log \sum_{R_{1,2}^{1}=u_M} 
W(\vrho^1,\vrho^2)
\exp \sqrt{N} \sum_{\ell\leq 2}y^{\ell}(\vrho^{\ell}).
\label{above}
\end{eqnarray}
If we define $\alpha=(\vrho^1,\vrho^2),$ define 
\begin{equation}
\A=\{(\vrho^1,\vrho^2): R_{1,2}^1=R(\vrho^1,\vrho^2) = u_M\},
\label{As}
\end{equation}
let
\begin{equation}
w_{\alpha}=W(\vrho^1,\vrho^2)/\sum_{R_{1,2}^1 = u_M} W(\vrho^1,\vrho^2),
\label{ws}
\end{equation}
and let
$z_i^{\ell}(\alpha) = z_i^{\ell}(\vrho^{\ell})$
and 
$y^{\ell}(\alpha) = y^{\ell}(\vrho^\ell),$
then (\ref{above}) can be rewritten as
\begin{eqnarray}
G_{M,N}
&=&
\frac{1}{N}
\e \log \sum_{\alpha\in\A} w_{\alpha}
\sum_{R_{1,2}^2 = u_N'}\exp
\sum_{\ell\leq 2}\sum_{i\leq N}\tau_i^{\ell}(z_i^{\ell}(\alpha) + h_{\ell})
\nonumber
\\
&-&
\frac{1}{N}\e \log \sum_{\alpha\in\A} 
w_{\alpha}
\exp \sqrt{N} \sum_{\ell\leq 2}y^{\ell}(\alpha).
\label{GM2}
\end{eqnarray}
Clearly, $G_{M,N}$ is written in the form of (\ref{GN}),
i.e. $G_{M,N} = G_N(u_N',\Omega_{\delta}),$
where the random overlap structure $\Omega_{\delta}$ is the collection
of (\ref{As}), (\ref{ws}), (\ref{zee}) and (\ref{ys}).
Equations (\ref{covzee}) and (\ref{covys}) imply that
the conditions (1) - (4) in the definition of ROSt
are satisfied with 
$
q_{\alpha,\beta}^{\ell,\ell'}=R(\vrho^{\ell},\bar{\vrho}^{\ell'})
$
where $\beta=(\bar{\vrho}^1,\bar{\vrho}^2).$
Since $q_{\alpha,\alpha}^{1,2}=R(\vrho^1,\vrho^2)=u_M$ for
$\alpha=(\vrho^1,\vrho^2)\in\A,$ we can take
$\delta=|u_M-u|$ which goes to $0$ as $M\to\infty.$
Equation (\ref{Plim}), therefore, implies
$$
\P(u)\geq \lim_{\delta\to 0 }\inf_{\Omega_{\delta}} 
G_N(u_N',\Omega_{\delta}).
$$
This finishes the proof of Theorem \ref{Th3}.
\qed

Lemma \ref{Lemma3} and Theorem \ref{Th3}, of course,
imply Theorem \ref{Th2}.

In  conclusion, we would like to note that the analogue of
the Aizenman-Sims-Starr variational principle is particularly
interesting because of the specific representations of the
random overlap structures (\ref{above}) or (\ref{GM2}).
We hope that the analysis of these structures will direct
toward what should be the Parisi ansatz for the limit $\P(u)$ in (\ref{limit})
or, at least, will provide some ideas in this direction.

\end{document}

%% file: invlat.tex
\font\tencmmib=cmmib10 \skewchar\tencmmib '60
\newfam\cmmibfam
\textfont\cmmibfam=\tencmmib

\def\bbox{\quad\hbox{\vrule \vbox{\hrule \vskip2pt \hbox{\hskip2pt
\vbox{\hsize=1pt}\hskip2pt} \vskip2pt\hrule}\vrule}}
\def\lessim{\ \lower4pt\hbox{$
\buildrel{\displaystyle <}\over\sim$}\ }
\def\gessim{\ \lower4pt\hbox{$\buildrel{\displaystyle >}
\over\sim$}\ }

\def\eps{\varepsilon}

\def\go0{\to 0}

\def\la{\langle}

\def\leftitem#1{\item{\hbox to\parindent{\enspace#1\hfill}}}

\def\qed{\hfill\break\rightline{$\bbox$}}

\def\ra{\rangle}

\def\sg{\sigma}

\def\sg2{\sigma^2}

\def\__{_{\infty}}